\documentclass[12pt]{article}
\usepackage{latexsym}
\usepackage{amsbsy}
\usepackage{graphics}
\usepackage{tikz}
\usepackage{pgf}
\usepackage{pgflibraryarrows}

\usepackage{amssymb}

\PassOptionsToPackage{hyphens}{url}
\usepackage{hyperref}

\usetikzlibrary{calc}   %  for co-ord calculations

\definecolor{greenish}{rgb}{0,0.5,0}
\definecolor{gold}{rgb}{0.85,.66,0}
\newtheorem{thm}{Theorem}
\newtheorem{cor}[thm]{Corollary}

\newcommand{\pf}{\noindent\textit{Proof.\ \ }}
\newcommand{\eopf}{\hfill\hspace*{4pt}\hfill
\fbox{\rule[-1pt]{0pt}{4pt}\hspace*{2pt}}}
\newcommand{\dd}{\ensuremath{.\,.}}

\newcommand{\vr}[2]{
\left( \begin{array}{c}#1\\#2\end{array} \right)}
\newcommand{\textvr}[2]{
\hbox{\footnotesize$\vr{#1}{#2}$}}

\newcommand{\minor}[1]{\,\|_{{}_{#1}}}
\newcommand{\mino}[1]{[#1]}    %  8/4/12

\newlength{\hatsize}
\newlength{\hataltitude}

\title{Binary functions, degeneracy, and alternating dimaps}   %  provisional title
\author{
G. E. Farr\thanks{Part of the work of this paper was presented at the Discrete Mathematics
Research Group, Monash University, June 2014.
This research was supported in part by ARC Discovery Grant DP110100957.}
\thanks{Email: \href{mailto:Graham.Farr@monash.edu}{\texttt{Graham.Farr@monash.edu}}}   \\
Faculty of Information Technology,   \\   %  *****
Monash University, Clayton, Victoria 3800,   \\
Australia  \\
}
\date{20 August 2017}     %  This file was started by taking the section on binary functions (used to be Section 4) out of my paper on `Minors for alternating dimaps' and starting a new paper based on it.  The source file for that other paper was started on 6 April 2012.  This file was created on 18 July 2016, with some revisions on 19 July 2016.  Further revisions were done 29 July 2017 - 5 Aug 2017 and 20 Aug 2017.
\begin{document}

\maketitle

\begin{abstract}
This paper continues the study of combinatorial properties of binary functions --- that
is, functions $f:2^E\rightarrow\mathbb{C}$ such that $f(\emptyset)=1$,
where $E$ is a finite set.  Binary functions have previously been shown to admit
families of transforms that generalise duality, including a trinity transform,
and families of associated minor operations that generalise deletion and
contraction, with both these families parameterised by the complex numbers.
Binary function representations exist for graphs (via the indicator functions of their
cutset spaces) and indeed arbitrary matroids (as shown by the author previously).
In this paper, we characterise degenerate elements --- analogues of loops and coloops ---
in binary functions, with respect to any pair of minor operations from our
complex-parameterised family.  We then apply this to study the relationship between
binary functions and Tutte's alternating dimaps, which also support a trinity transform
and three associated minor operations.  It is shown that only the simplest alternating
dimaps have binary representations of the form we consider, which seems to be the
most direct type of representation.  The question of whether there exist other,
more sophisticated types of binary function representations for alternating dimaps is
left open.
\end{abstract}

\textit{2010 Mathematics Subject Classification:}
Primary: 05B99, 05C83;
Secondary: 05B35, 05C10, 05C20, 05C31, 05C50, 05C99, 05E99.   \\   \\   \\

\section{Introduction}

Duality is a pervasive theme in discrete mathematics.  It runs strongly through
planar graph theory, matroid theory, mathematical programming, and (via the
Hadamard transform) information and coding theory.  In some of these contexts,
there are also \textit{minor} operations (local ``reductions'' based on a specific
element) that are dual to each other --- the
main example being deletion and contraction for graphs and matroids.  These
operations are important in characterising structures with various properties and
in the theory of enumeration.

In some combinatorial systems, transforms of higher order than duality ---
\textit{triality} or \textit{trinity}, which have order three rather than two --- are important.
These may also have minor operations of some type, now three in number.
The main object types with a trinity transform and associated minors
known to the author are alternating dimaps (which go back to
Tutte \cite{tutte48}, with minor operations in \cite{farrXX}), binary functions \cite{farr2013},
multimatroids (including isotropic systems) \cite{bouchet87,bouchet98} and the related
transition matroids \cite{traldi2015}.  This paper continues our ongoing study of
combinatorial trinities and associated minor operations \cite{farr2013,farrXX},
by exploring the relationship between alternating dimaps and binary functions.

An \textit{alternating dimap} is an
orientably embedded directed graph in which, around each vertex, the incident edges
are directed alternately into, and out of, the vertex.  Alternating dimaps were introduced
by Tutte \cite{tutte48} as part of his work on dissecting equilateral triangles into equilateral
triangles.

A \textit{binary function} is a complex-valued function defined on all subsets of a set that
takes value 1 on the empty set.  The prototypical example is the indicator function of the
cutset space of a graph, or of a binary linear space (i.e., a cocircuit space of a binary
matroid).  Other examples come from indicator functions of
powerful sets  \cite{farr-wang2017} or indeed any collection of finite
sets that includes the empty set.  In these cases, the binary functions
only need to be $\{0,1\}$-valued.  But larger ranges allow binary functions to
generalise \textit{all} matroids \cite{farr93},
to support minor operations more general than just deletion and contraction \cite{farr04},
and to support transforms other than duality, such as trinity
and higher order transforms \cite{farr2013}.

One reason for our interest in binary functions is that they have their own
Tutte-Whitney functions (not always strictly polynomials, since the exponents of the
variables may not necessarily be integers).
These were introduced in \cite{farr93}, and shown there to contain
functions of independent interest such as the weight enumerator of an arbitrary
(not necessarily linear) code,
Oxley and Welsh's clutter reliability (or percolation probability) \cite{oxley-welsh79a},
and Kung's generalised chromatic polynomial \cite{kung80}.
A natural generalisation of the Potts model partition function to binary functions
was similarly treated in \cite{farr04}.
The Tutte-Whitney polynomials (or functions)
of binary functions and their duals (i.e., Hadamard
transforms) were found in \cite{farr04} to be just two members of a whole family of generalised Tutte-Whitney functions, and these were found in \cite{farr07b} to contain
the partition function of the symmetric Ashkin-Teller model, which cannot be found from
the usual Tutte-Whitney polynomials.

To better understand the combinatorics of binary functions, we need to understand
their degenerate elements.  An element is \textit{degenerate}
if all possible minor operations on it give the same result.
For example, the degenerate elements
of a graph or a matroid are its loops and coloops.  These play a fundamental role:
they constitute the singleton components of the matroid;
they may be regarded as base cases
for the recursive definitions of Tutte-Whitney polynomials; and they must be treated as
special cases in numerous proofs.  Degenerate elements of alternating
dimaps are triloops, in the terminology of \cite{farrXX}.
We characterise degenerate elements for binary functions
in \S\ref{sec:bin-fins-degen}.  We then apply the characterisation
to study the relationship between binary functions and alternating dimaps.

Some alternating dimaps certainly cannot be represented by binary functions,
or indeed any of the other types of objects mentioned in the second paragraph,
since alternating dimap minor operations may not commute, unlike those in the
other settings.
In \S\ref{sec:strict-bin-reprns}
we determine those alternating dimaps that can be represented faithfully by binary functions.

\section{Alternating dimaps}

An \textit{alternating dimap} is a directed graph
without isolated vertices whose components are each
2-cell-embedded in a separate orientable 2-manifold, such that
for each vertex $v$,
the edges incident with $v$ are
directed alternately into, and out of, $v$ (going around $v$ in the embedding).
So each vertex in an alternating dimap has equal indegree and outdegree.
Loops and/or multiple edges are allowed, but coloops are not possible.
The \textit{empty alternating dimap} has no
vertices, edges or faces.  Alternating dimaps were introduced by Tutte and his
collaborators, who developed their theory
in \cite{tutte48,brooks-etal75,tutte75,tutte99}; see also \cite[\S10.3]{brooks-etal40} and \cite[Ch.\ 4]{tutte98}.  Work on the topic by others includes \cite{berman80,cori-penaud80,kalman2013}, with much of this work focusing on Tutte's Tree Trinity Theorem.
For a brief summary of the main elements of the theory that concern us here,
and for a review of the history and related work, see \cite{farrXX}.

If $G$ is an alternating dimap then $kG$ is the disjoint union
of $k$ copies of $G$.  Each of these copies is regarded as being
embedded in separate surfaces.

%An edge $e$ from $u$ to $v$ is sometimes written $e(u,v)$.

A face is \textit{clockwise} or
\textit{anticlockwise} according to the direction of the edges
around it.  (All faces are of one of these two types).
If two faces share a common edge, then one of the faces
is clockwise and the other is anticlockwise.
The \textit{left successor} (respectively, \textit{right successor})
of an edge $e$ is the next edge
after $e$, going around its anticlockwise (resp., clockwise) face
in the direction given by $e$ (i.e., anticlockwise or clockwise, respectively).

If $G$ is an alternating dimap, then its \textit{trial} $G^{\omega}$, introduced
by Tutte \cite{tutte48}, may be defined as
follows.  The vertices of $G^{\omega}$ correspond to the clockwise faces of $G$;
write $C_u$ for the clockwise face in $G$ represented by vertex $u$ in $G^{\omega}$.
There is a directed edge $(u,v)$ from $u$ to $v$ in $G^{\omega}$ whenever
there is a vertex $a$ in $G$ belonging to both $C_u$ and $C_v$ such that $C_u$ is
the next clockwise face after $C_v$, going clockwise around $a$.  This edge is defined
to be the image $e^{\omega}$, under triality, of the edge $e$ of $C_v$ that goes into $a$
and whose left successor is an edge of $C_u$ going out of $a$.  See \cite{farrXX}
for a more detailed treatment.

Tutte showed that $((G^{\omega})^{\omega})^{\omega}=G$.
The symbol $\omega$ is taken to satisfy $\omega^3=1$, so we can write, for example,
$((G^{\omega})^{\omega})^{\omega}=(G^{\omega^2})^{\omega}=G^{\omega^3}=G^1=G$
At times it is natural to put $\omega=\exp(2\pi i/3)$.

In \cite{farrXX}, the author introduced three \textit{reductions}, or \textit{minor operations},
that may be done to any edge $e$ in an alternating dimap,
and which are the analogues in this
context of deletion and contraction in graphs.  More specifically, they are the analogues
of the \textit{surface minor} versions of deletion and contraction, which apply to embedded
graphs rather than abstract graphs.  The first of these three reductions
is \textit{1-reduction} or
\textit{contraction}, which behaves exactly like the standard contraction operation on
embedded graphs.  In particular, if $e$ is a loop at a vertex $v$ such that $e$
does not constitute a face in its
own right, then contraction of $e$ causes its incident vertex to split into two copies of
itself, one for each of the two sides into which $e$ divides the neighbourhood of $v$
in the surface.  This split will either increase the number of components or reduce the
genus.  The second reduction is $\omega$\textit{-reduction}, in which the left successor
of $e$ is changed so that it starts at the tail of $e$ (and its head is unchanged) and $e$
is deleted.  The third reduction is $\omega^2$ reduction, which is defined as for
$\omega$-reductions except using the right successor rather than the left successor.
A \textit{minor} of an alternating dimap is another alternating dimap obtained from the first
by some sequence of reductions.

The relationship between triality and minors is described by the following result.

\begin{thm}
\cite[Theorem 2.2]{farrXX}   %  was Thm 4 in preprint
\label{thm:trial-minor}
If $e\in E(G)$ and $\mu,\nu\in\{1,\omega,\omega^2\}$ then
\[
G^{\mu}\mino{\nu}e^{\mu} = (G\mino{\mu\nu}e)^{\mu} .
\]
\end{thm}

This extends the classical relationship between duality, deletion and contraction,
under which $G^*\setminus e  = (G/e)^*$ and $G^*/e = (G\setminus e)^*$.

The relationship between triality and minor operations for
alternating dimaps given by Theorem \ref{thm:trial-minor}
is reminiscent of properties of binary
functions found by the author in \cite{farr2013}.  
One outcome of the present paper is to establish how close this connection is.
We will determine those alternating dimaps which can be
represented, in a certain faithful manner, by binary functions.
It is found that alternating dimaps and binary functions actually do not
have much in common.

A special role in the theory of alternating dimaps is played by various types of special
edges, which are analogous to loops and coloops in graphs but more diverse.
An \textit{ultraloop} is a loop forming a component in its own right.
A \textit{1-loop} is an edge whose head has indegree = outdegree = 1 (and which need
not necessarily be a loop in the conventional sense).
An $\omega$-loop is a loop forming an anticlockwise face of size 1, while
an $\omega^2$-loop is a loop forming a clockwise face of size 1.
An ultraloop is therefore also a 1-loop, an $\omega$-loop, and an $\omega^2$-loop.
A \textit{triloop} is an edge that is a 1-loop, an $\omega$-loop or an $\omega^2$-loop,
and it is \textit{proper} if it is not also an ultraloop.
The triloops are precisely those edges such that the three reductions on it all give
identical minors.

For abstract graphs, loops and coloops are the only types of special edges.  But
embedded graphs have other types of special edges too, which are not necessarily
degenerate:
a \textit{semiloop} (respectively, a \textit{semicoloop}) is either a loop
(resp., coloop) or an edge whose
contraction (resp., deletion) either increases the number of components or decreases
the genus.  These too have analogues for alternating dimaps.

A \textit{1-semiloop} is just an ordinary loop.
An $\omega$-semiloop (respectively, an $\omega^2$-semiloop)
is an $\omega^2$-loop (resp., an $\omega$-loop) or
an edge for which $\omega^2$-reduction (resp., an $\omega$-reduction)
either increases the
number of components or reduces the genus.
A semiloop is \textit{proper} if it is not also a triloop.

\section{Binary functions}
\label{sec:bin-fns}

In this
section we briefly summarise some of the theory of binary functions
developed by the author in \cite{farr93,farr04,farr07a,farr07b,farr2013}.
We restrict attention to aspects that are relevant to this paper, and so
focus on \cite{farr2013}.

Let $E$ be a finite set, with $m=|E|$.
A \textit{binary function} with \textit{ground set} $E$
and \textit{dimension} $m$
is a function $f:2^E\rightarrow\mathbb{C}$
such that $f(\emptyset)=1$.  Equivalently, we regard it as
a $2^{m}$-element complex vector $\mathbf{f}$
whose elements are indexed
by the subsets of $E$ and whose first element (indexed by $\emptyset$) is 1.
(The restriction $f(\emptyset)=1$ was not imposed as part of the definition
in earlier work \cite{farr93,farr04,farr07a,farr07b}, but all scalar
multiples of a binary function are equivalent for our purposes, and
we have always been most interested in the cases where $f(\emptyset)\not=0$.)
If $f$ is a binary function then $E(f)$ denotes its ground set.

We often represent a subset $X\subseteq E$ by its characteristic
vector $\mathbf{x}\in\{0,1\}^E$, with $x_e=1$ if $e\in X$ and $x_e=0$ otherwise;
here, the set $E$ indexes the positions in the characteristic vector.
Since $\mathbf{x}$ may be thought of as a binary string, it may also
be taken to be the binary representation of a number $x$ such that
$0\le x\le 2^{m}-1$, using $E=\{0,\ldots,m-1\}$.
These numbers give the order in which the subsets of $E$ are listed, as indices
of the entries of the vector $\mathbf{f}$.
With this notation, $f(X)$ may also be written
$f_{\mathbf{x}}$ or $f_x$.
In particular, $f(\emptyset)=f_{(0,\ldots,0)}=f_0=1$.
We write $0_k$ for the sequence of $k$ 0s,
and sometimes drop the subscript $k$ when it is clear from the context.

The definition was motivated by indicator functions of linear spaces
over GF(2), especially of cutset spaces of graphs: if $N$ is a matrix
over GF(2) whose columns are indexed by $E$ (such as the incidence
matrix of a graph, or the matrix representation of a binary matroid),
then the indicator function of the rowspace of $N$ takes value 1
on a set $X\subseteq E$ if the characteristic vector of $X$ belongs
to the rowspace of $N$, and takes value 0 otherwise.

If $f, g:2^E\rightarrow\mathbb{C}$ and there exists a constant
$c\in\mathbb{C}\setminus\{0\}$ such that $f(X) = c g(X)$
for all $X\subseteq E$, then we write $f\simeq g$.

Define
\[
M(\mu) := \frac{1}{2\sqrt{2}}
\left(
\begin{array}{cc}
\sqrt{2}+1+(\sqrt{2}-1)\mu  &  1-\mu   \\
1-\mu  &  \sqrt{2}-1+(\sqrt{2}+1)\mu
\end{array}
\right) .
\]
The $\mu$\textit{-transform} of $\mathbf{f}$, denoted by
$L^{[\mu]}\mathbf{f}$, is given by
\[
L^{[\mu]}\mathbf{f} := M(\mu)^{\otimes m}\mathbf{f} ,
\]
where the $2^m\times2^m$ matrix on the right is the $m$-th Kronecker
power of $M(\mu)$.

When $\mu=1$, we have the identity transform, while when $\mu=-1$,
we have a scalar multiple of the Hadamard transform.  It is well known that
the Hadamard transform takes the indicator function of a linear space
to a scalar multiple of the indicator function of its dual, from which
it follows that the indicator functions of the cutset and circuit spaces
of a graph are related by the Hadamard transform in the same way.
It was shown in \cite{farr93} that general matroid duality is also
described by the Hadamard transform.

It is easy to show that $M(\mu_1\mu_2)=M(\mu_1)M(\mu_2)$ for all $\mu_1,\mu_2$,
see \cite{farr2013}.
It follows (using the mixed-product property for the Kronecker product)
that composition of the $L^{[\mu]}$ transforms corresponds to multiplication
of their $\mu$ parameters: $L^{[\mu_1]} L^{[\mu_2]} = L^{[\mu_1\mu_2]}$,
from \cite[Theorem 2]{farr2013}.  At this point, readers may ask: what
happens when $\mu=\omega$?  We look at this shortly.

Suppose $E=\{e_0,\ldots,e_{m-1}\}$.

Let $b\in\{0,1\}$.  We use $\mathbf{f}_{b\bullet}$ as shorthand for the
vector of length $2^{m-1}$, with elements indexed by
subsets of $E\setminus e_0$, whose $X$-element
is $f(X)$, if $b=0$, or $f(X\cup\{e_0\})$, if $b=1$
(for $X\subseteq E\setminus\{e_0\}$).
We define $\mathbf{f}_{\bullet b}$ in the same way,
except that we use $e_{m-1}$ instead of $e_0$ throughout.
The vectors $\mathbf{f}_{0\bullet}$ and $\mathbf{f}_{1\bullet}$
give the top and bottom halves, respectively, of $\mathbf{f}$,
while $\mathbf{f}_{\bullet 0}$ and $\mathbf{f}_{\bullet 1}$
give the elements in even and odd positions, respectively,
of $\mathbf{f}$.

Let $I_l$ denote
the $l\times l$ identity matrix.
If $e\in E$, then the $[\mu]$\textit{-minor} of $\mathbf{f}$ by $e$ is
the $2^{m-1}$-element
vector $\mathbf{f}\minor{[\mu]}e$, with entries indexed by
subsets of $E\setminus\{e\}$, given by
\begin{equation}
\label{eq:bin-fn-minor}
\mathbf{f}\minor{[\mu]}e_i ~ := ~ c \cdot  ~(~ I_2^{\otimes i} \otimes
( ~~ 1 ~~~ \frac{1+\mu}{\sqrt{2}+1-(\sqrt{2}-1)\mu} ~~ )
\otimes I_2^{\otimes(m-i-1)} ~)~ \mathbf{f} ,
\end{equation}
where $c$ is such that the $\emptyset$-element
of $\mathbf{f}\minor{[\mu]}e_i$ is 1.

%Put $(\mu_0,\mu_1,\mu_2)=(1,\omega,\omega^2)$ and, for each $j\in\{0,1,2\}$,
For any $\mu$, define $\lambda=\lambda(\mu)$ by
\begin{equation}
\label{eq:lambda-mu}
\lambda(\mu) := \frac{1+\mu}{\sqrt{2}+1-(\sqrt{2}-1)\mu} .
\end{equation}
Then $\mathbf{f}\minor{[\mu]}e_i$ is a scalar multiple of
$(~ I_2^{\otimes i} \otimes
( ~~ 1 ~~~ \lambda ~~ )
\otimes I_2^{\otimes(m-i-1)} ~)~ \mathbf{f}$.

When $f$ is the indicator function of the cutset space of a graph,
the minor $f\minor{[\mu]}e$ amounts to deletion when $\mu=1$ and
contraction when $\mu=-1$.  See \cite[\S2,\S6]{farr2013}, and also
\cite{farr04} for the first definition of generalised minor operations
interpolating between deletion and contraction (albeit with a different
parameterisation to that used here and in \cite{farr2013}).
This work has its roots in \cite{farr93}, where deletion and contraction
are expressed in terms of indicator
functions of cutset spaces, and these operations are extended to
general binary functions.

The relationship between transforms and minors for binary functions 
is as follows.

\begin{thm}
\cite[Theorem 6.1]{farr2013}
\label{thm:transform-minor-bin-fn}
If $e\in E(f)$ and $\mu,\nu\in\mathbb{C}$ then
\[
(L^{[\mu]}f)\minor{[\nu]}e  \simeq L^{[\mu]}(f\minor{[\mu\nu]}e) .
\]
\end{thm}
This may be compared with Theorem \ref{thm:trial-minor}.
In particular, we have
\begin{eqnarray*}
(L^{[\omega]}f)\minor{[1]}e  &  \simeq  &
L^{[\omega]}(f\minor{[\omega]}e) ,   \\
(L^{[\omega]}f)\minor{[\omega]}e  &  \simeq  &
L^{[\omega]}(f\minor{[\omega^2]}e) ,   \\
(L^{[\omega]}f)\minor{[\omega^2]}e  &  \simeq  &
L^{[\omega]}(f\minor{[1]}e) .
\end{eqnarray*}

\section{Alternating dimaps and binary functions}

The relationship described above between the transform $L^{[\omega]}$
(called the \textit{trinity transform} \cite{farr2013} or
\textit{triality transform}) and
the minor operations 
for binary functions follows the same pattern as the relationships
between triality and minors for alternating dimaps,
given in Theorem \ref{thm:trial-minor}.
It is natural to ask what connection there may be between the two.

For binary functions, the minor operations always
commute \cite[Lemma 4]{farr04}.  In fact, that result implies that
every binary function is \textit{totally reduction-commutative}, meaning that
for any set of reductions $\bullet\minor{[\mu_i]}e_i$ on distinct edges $e_i$,
any ordering of the reductions
gives the same minor (borrowing some terminology from \cite{farrXX}).
But, as we saw in \cite[\S3]{farrXX},
the minor operations for alternating dimaps do not always commute.
It follows that alternating dimaps, along with triality and minor
operations, cannot be represented faithfully by binary functions with
their trinity transform and minor operations described above. 

Nonetheless, we can ask if there is a subclass of alternating dimaps
which can be represented faithfully by binary functions in this way.
For this to occur, this subclass must consist only of alternating
dimaps that are totally reduction-commutative.  Such alternating
dimaps were characterised in \cite[Theorem 17]{farrXX}; the subclass
we seek must be a subset of those.

Later we will give a definition of faithful representation by
binary functions, and determine when such a representation is possible.
To do the latter, it will help to characterise those binary functions
for which any $\mu$-reduction ($\mu\in\{1,\omega,\omega^2\}$),
on any element of the ground set, gives the same result.

\section{Degeneracy for binary functions}
\label{sec:bin-fins-degen}

We now extend the term ``degenerate'' to any type of
combinatorial structure on which some kind of minor operations are defined: an element
is \textit{degenerate} if all reductions of it, using minor operations, give the same object.
Any real understanding of a particular type of combinatorial structure with minors can be
expected to depend, in part, on understanding the degenerate elements.

For alternating dimaps, the degenerate elements are the triloops (including ultraloops).

In this section, we consider degenerate elements for binary functions under the three
minor operations.  To do this, we need some more notation.

Throughout, we write
\[
\mathbf{i}=\vr{1}{0}, ~~~~
\mathbf{j}=\vr{0}{1}, ~~~~
H = (\mathbf{h}_0,\ldots,\mathbf{h}_{k-1})
\in \{\mathbf{i},\mathbf{j}\}^{\{0,\ldots,k-1\}} .
\]
For each $i$,
\[
H^{(i)} = (\mathbf{h}_0,\ldots,\mathbf{h}_{i-1},
\mathbf{h}_{i+1},\ldots,\mathbf{h}_{k-1})
\]
is the sequence obtained from $H$ by omitting the term indexed by $i$.

For each $H$, define the sequence $G=G(H)=(g_0,\ldots,g_{k-1})$ by
\[
g_i = \left\{
\begin{array}{ll}
0,  &  \hbox{if $\mathbf{h}_i = \mathbf{i}$};   \\
1,  &  \hbox{if $\mathbf{h}_i = \mathbf{j}$}.
\end{array}
\right.
\]
The sequence obtained from this by omitting the term indexed by $i$ is
\[
G^{(i)} = (g_0,\ldots,g_{i-1},g_{i+1},\ldots,g_{k-1}) .
\]
The subsequence $(g_{i_1},\ldots,g_{i_2})$ of $G$ is denoted by
$G[i_1\dd i_2]$.

If $b\in\{0,1\}$, then $G:i\leftarrow b$ denotes the sequence obtained by inserting $b$
between the $i$-th and $(i+1)$-th elements of $G$:
\[
G:i\leftarrow b ~=~ (g_0, \ldots, g_{i-1}, b, g_i, \ldots, g_{k-1}) .
\]
The two-element vector $\mathbf{f}_{G:i}$ is defined by
\[
\mathbf{f}_{G:i} = \vr{f_{G:i\leftarrow0}}{f_{G:i\leftarrow1}} .
\]

Write $\mathbf{u}$ for a $2^k$-element vector indexed by the numbers
$0,\ldots,2^{k-1}$ --- or, equivalently, by vectors of $k$ bits, or by subsets of
$\{0,\ldots,k-1\}$.

For a given $G$, we write $u_{G}$ for the entry of $\mathbf{u}$ whose
index has binary representation given by $G$, i.e., whose index is
$\sum_{i=0}^{k-1}g_i2^{k-1-i}$.
%  Write $\mathbf{u}_{G}^{(i)}$ for the vector
%  \[
%  \vr{u_{G[0\dd i-1],0,G[i+1\dd k-1]}}{u_{G[0\dd i-1],1,G[i+1\dd k-1]}} ,
%  \]
%  where each subscript is a sequence of $k$ bits, interpreted for indexing
%  purposes as a binary number.

It is routine to show that, if $m\ge1$ and
$u$ is a (vector representation of a)
binary function with ground set of size $m-1$, then
\begin{equation}
\label{eq:u-sum}
\mathbf{u} = \sum_{H}
(~ \mathbf{h}_0 \otimes \cdots \otimes
\mathbf{h}_{m-2} ~) ~ u_{G(H)} .
\end{equation}
If $m=1$ then there is a single $H$ to sum over, consisting of the
empty sequence, and the empty product
$\mathbf{h}_0 \otimes \cdots \otimes
\mathbf{h}_{m-2}$ is the trivial single-element vector $~(~ 1 ~)$.
Also $G=G(H)$ is the empty bit-sequence, representing the number 0, and
$u_G=u_0=1$, so $\mathbf{u} = (~1~)$, as expected.

%  Also, though I'm not sure yet if I'll need it, we have ********
%  \[
%  \mathbf{u} = \sum_{H^{(i)}}
%  \mathbf{h}_0 \otimes \cdots \otimes
%  \mathbf{h}_{i-1} \otimes \mathbf{u}_{G}^{(i)} \otimes
%  \mathbf{h}_{i+1} \otimes \cdots \otimes
%  \mathbf{h}_{m-2} .
%  \]

\begin{thm}
\label{thm:reduce-f-to-u-soln}
Let $\mu_1,\mu_2\in\mathbb{C}\setminus\{3+2\sqrt{2}\}$ be distinct.
Suppose $\mathbf{f}$ and $\mathbf{u}$ are binary functions of dimension $m$
and $m-1$ respectively.  Then
\[
\mathbf{f}\minor{[\mu_1]}e_i = \mathbf{f}\minor{[\mu_2]}e_i = \mathbf{u}
\]
if and only if
for all $G\in\{0,1\}^{\{0,\ldots,m-2\}}$ and all $b\in\{0,1\}$,
\begin{equation}
\label{eq:reduce-f-to-u-soln}
f_{G:i\leftarrow b}  =
f_{0:i\leftarrow b} \, u_{G} \,.
\end{equation}
\end{thm}

\pf
Let us write the hypothesis as a set of equations, using
(\ref{eq:bin-fn-minor}) and (\ref{eq:lambda-mu}).  The condition that
$\mathbf{f}\minor{[\mu_j]}e_i=\mathbf{u}$ for all $j\in\{1,2\}$ is equivalent to
the assertion that, for each such $j$, there exists $c_{ij}$ such that
\begin{equation}
\label{eq:reduce-f-to-u}
(~ I_2^{\otimes i} \otimes
( ~~ 1 ~~~ \lambda_j ~~ )
\otimes I_2^{\otimes(m-i-1)} ~)~ \mathbf{f} = c_{ij} \, \mathbf{u} \,,
\end{equation}
where $\lambda_j := \lambda(\mu_j)$.  Note that $\lambda_1\not=\lambda_2$.

Put
\[
R = \left(
\begin{array}{cc}
1 & \lambda_1   \\
1 & \lambda_2
\end{array}
\right)
~~~~~~ \hbox{and} ~~~~~~
\mathbf{c}_i =
\left(
\begin{array}{c}
c_{i1}  \\  c_{i2}
\end{array}
\right) .
\]
The equations (\ref{eq:reduce-f-to-u}) may be written (using (\ref{eq:u-sum})),
\begin{equation}
\label{eq:reduce-f-to-u-new}
(~ I_2^{\otimes i} \otimes  R
\otimes I_2^{\otimes(m-i-1)} ~)~ \mathbf{f}
= \sum_{H}
(~ \mathbf{h}_0 \otimes \cdots \otimes \mathbf{h}_{i-1} \otimes
\mathbf{c}_i \otimes
\mathbf{h}_{i} \otimes \cdots \otimes \mathbf{h}_{m-2} ~) ~ u_{G} \,.
\end{equation}
Here and below we write $G=G(H)$ for brevity.
We may write
\[
\mathbf{f} =
\sum_{H}  (~ \mathbf{h}_0 \otimes \cdots \otimes \mathbf{h}_{i-1} \otimes
\mathbf{f}_{G:i} \otimes
\mathbf{h}_{i} \otimes \cdots \otimes \mathbf{h}_{m-2} ~) ,
\]
so the left-hand side of (\ref{eq:reduce-f-to-u-new}) is
\begin{eqnarray*}
\lefteqn{
(~ I_2^{\otimes i} \otimes  R \otimes I_2^{\otimes(m-i-1)} ~)~
\sum_{H}  (~ ( \bigotimes_{k=0}^{i-1} \mathbf{h}_k ) \otimes
\mathbf{f}_{G:i} \otimes
(\bigotimes_{k=i}^{m-2} \mathbf{h}_{k} ) ~) }   \\
& = &
\sum_{H}  (~ I_2^{\otimes i} \otimes  R \otimes I_2^{\otimes(m-i-1)} ~)~
(~ ( \bigotimes_{k=0}^{i-1} \mathbf{h}_k ) \otimes
\mathbf{f}_{G:i} \otimes
(\bigotimes_{k=i}^{m-2} \mathbf{h}_{k} ) ~)   \\
& = &
\sum_{H}  
(~ ( \bigotimes_{k=0}^{i-1} I_2\,\mathbf{h}_k ) \otimes
R\,\mathbf{f}_{G:i} \otimes
(\bigotimes_{k=i}^{m-2} I_2\,\mathbf{h}_{k} ) ~)   \\
& = &
\sum_{H}  
(~ ( \bigotimes_{k=0}^{i-1} \mathbf{h}_k ) \otimes
R\,\mathbf{f}_{G:i} \otimes
(\bigotimes_{k=i}^{m-2} \mathbf{h}_{k} ) ~) .
\end{eqnarray*}
Setting this equal to the right-hand side of (\ref{eq:reduce-f-to-u-new}) and equating
appropriate elements, we find that (\ref{eq:reduce-f-to-u-new}) is equivalent to
\begin{equation}
\label{eq:reduce-f-to-u-soln1}
R \, \mathbf{f}_{G:i}  =  \mathbf{c}_i \, u_G  \,, ~~~~
\hbox{for all $G\in\{0,1\}^{\{0,\ldots,m-2\}}$}.
\end{equation}
%for all $G\in\{0,1\}^{\{0,\ldots,m-2\}}$.

When $G=0$, (\ref{eq:reduce-f-to-u-soln1}) and $u_0=1$ give
\begin{equation}
\label{eq:Rf0ici}
R \, \mathbf{f}_{0:i}  =  \mathbf{c}_i .
\end{equation}

For all $G\in\{0,1\}^{\{0,\ldots,m-2\}}$, we have
\begin{eqnarray*}
R \, \mathbf{f}_{G:i}  =  \mathbf{c}_i \, u_G
& \Longleftrightarrow &
R \, \mathbf{f}_{G:i} = R \, \mathbf{f}_{0:i} \, u_G
~~~~~ \hbox{(using (\ref{eq:Rf0ici}))}   \\
& \Longleftrightarrow &
\mathbf{f}_{G:i} = \mathbf{f}_{0:i} \, u_G
~~~~~ \hbox{(since $R$ is invertible, because the $\lambda_j$ are distinct)}   \\
& \Longleftrightarrow &
\forall b\in\{0,1\}: ~~
f_{G:i\leftarrow b}  =
f_{0:i\leftarrow b} \, u_{G} \,.
\end{eqnarray*}
%Using this with (\ref{eq:reduce-f-to-u-soln1}) gives, for any $G$,
%\[
%R \, \mathbf{f}_{G:i} = R \, \mathbf{f}_{0:i} \, u_G .
%\]
%Since $R$ is invertible (because the $\lambda_j$ are distinct),
%\[
%\mathbf{f}_{G:i} = \mathbf{f}_{0:i} \, u_G .
%\]
%Hence, for all $b\in\{0,1\}$ and all $G$, 
%\[
%f_{G:i\leftarrow b}  =
%f_{0:i\leftarrow b} \, u_{G} \,.
%\]
%
%Conversely, if this last equation holds for all $b$ and $G$, we obtain
%\[
%\mathbf{f}_{G:i} = \mathbf{f}_{0:i} \, u_G .
%\]
%and therefore
%\[
%R \, \mathbf{f}_{G:i} = R \, \mathbf{f}_{0:i} \, u_G .
%\]
%Set $\mathbf{c}_i  :=  R \, \mathbf{f}_{0:i}$,
%so that, for all $G$, we have
%\[
%R \, \mathbf{f}_{G:i}  =  \mathbf{c}_i \, u_G  \,.
%\]
%This is just (\ref{eq:reduce-f-to-u-soln1}), which we showed to be equivalent to 
%(\ref{eq:reduce-f-to-u-new}) and hence to (\ref{eq:reduce-f-to-u}).
\eopf   \\

Using this result, we can characterise degenerate elements in a binary function
as those elements $e_i$ such that the binary function entries whose index vectors
differ only in position $i$ have constant ratio or are both zero.

\begin{cor}
Let $\mu_1,\mu_2\in\mathbb{C}\setminus\{3+2\sqrt{2}\}$ be distinct.
Let $\mathbf{f}$ be a binary function of dimension $m$.  The following are equivalent.
\begin{enumerate}
\item[(a)]
The element $e_i$ is degenerate in $\mathbf{f}$ with respect to the minor operations
$\bullet\minor{[\mu_1]}$ and $\bullet\minor{[\mu_2]}$.
\item[(b)]
For all $G\in\{0,1\}^{\{0,\ldots,m-2\}}$
\begin{equation}
\label{eq:degeneracy-products}
f_{G:i\leftarrow 1} f_{0:i\leftarrow 0}  =
f_{G:i\leftarrow 0} f_{0:i\leftarrow 1} .
\end{equation}
\item[(c)]
For all $G\in\{0,1\}^{\{0,\ldots,m-2\}}$ either
\begin{equation}
\label{eq:degeneracy-ratios}
\frac{f_{G:i\leftarrow 1}}{f_{G:i\leftarrow 0}}  =
\frac{f_{0:i\leftarrow 1}}{f_{0:i\leftarrow 0}}
\end{equation}
or
\begin{equation}
\label{eq:degeneracy-zeros}
f_{G:i\leftarrow 0} = f_{G:i\leftarrow 1}  = 0 .
\end{equation}
\end{enumerate}
\end{cor}

\pf
((a) $\Leftrightarrow$ (b)):
Degeneracy means that
\[
\mathbf{f}\minor{[\mu_1]}e_i = \mathbf{f}\minor{[\mu_2]}e_i =: \mathbf{u} ,
\]
where $\mathbf{u}$ has dimension $m-1$.
By Theorem \ref{thm:reduce-f-to-u-soln}, degeneracy is equivalent to the assertion that,
for all $G\in\{0,1\}^{\{0,\ldots,m-2\}}$ and all $b\in\{0,1\}$,
\begin{equation}
\label{eq:fGf0uG}
f_{G:i\leftarrow b}  =  f_{0:i\leftarrow b} \, u_{G} \,.
\end{equation}
Putting $b=0$ we have $u_G = f_{G:i\leftarrow 0} / f_{0:i\leftarrow 0}$;
substituting this into (\ref{eq:fGf0uG}) with $b=1$ gives (\ref{eq:degeneracy-products}).
Conversely, if (\ref{eq:degeneracy-products}) holds for all $G\in\{0,1\}^{\{0,\ldots,m-2\}}$,
we put $u_G := f_{G:i\leftarrow 0} / f_{0:i\leftarrow 0}$ for each $G$, which gives
$f_{G:i\leftarrow 0} = f_{0:i\leftarrow 0} u_G$ immediately, and
$f_{G:i\leftarrow 1} = f_{0:i\leftarrow 1} u_G$ after substitution
into (\ref{eq:degeneracy-products}).
Since (\ref{eq:fGf0uG}) then holds for all $G$ and $b$,
we conclude that $e_i$ is degenerate, by Theorem \ref{thm:reduce-f-to-u-soln}.

(c) $\Rightarrow$ (b) is immediate.

((b) $\Rightarrow$ (c)):
Suppose (\ref{eq:degeneracy-products}) holds for all $G$.
Then $f_{G:i\leftarrow 1} = f_{G:i\leftarrow 0} = 0$ or
$f_{G:i\leftarrow 1} = f_{0:i\leftarrow 1} = 0$ or both sides of
(\ref{eq:degeneracy-products}) are nonzero, in which case we have
\[
\frac{f_{G:i\leftarrow 1}}{f_{G:i\leftarrow 0}}  =
\frac{f_{0:i\leftarrow 1}}{f_{0:i\leftarrow 0}} .
\]
But this last equation subsumes $f_{G:i\leftarrow 1} = f_{0:i\leftarrow 1} = 0$, so
we have (c).
\eopf   \\

%
%f_{G:i\leftarrow 1} f_{0:i\leftarrow 0}  = 
%f_{G:i\leftarrow 0} f_{0:i\leftarrow 1} .
%
%
%If $u_G=0$ then $f_{G:i\leftarrow 0} = f_{G:i\leftarrow 1} = 0$,
%establishing (\ref{eq:degeneracy-zeros}).  So suppose $u_G\not=0$.
%We also have $f_{0:i\leftarrow 0}=1$, since $\mathbf{f}$ is a binary function
%(so $f(\emptyset)=1$).  Therefore $f_{G:i\leftarrow 0}=u_{G}\not=0$.  Using
%(\ref{eq:fGf0uG}) again with $b=1$, we have
%\[
%f_{G:i\leftarrow 1}
%=  f_{0:i\leftarrow 1} \, u_{G}
%=  f_{0:i\leftarrow 1} \, f_{G:i\leftarrow 0} ,
%\]
%giving (\ref{eq:degeneracy-ratios}),
%using $f_{0:i\leftarrow 0}=1$ and $f_{G:i\leftarrow 0}\not=0$.
%
%Now suppose that either (\ref{eq:degeneracy-ratios}) or (\ref{eq:degeneracy-zeros})
%holds, for all $G$.  Consider any particular $G$.
%Put $u_G=f_{G:i\leftarrow 0}$.
%If (\ref{eq:degeneracy-ratios}) holds for $G$, then
%\[
%f_{G:i\leftarrow 1}  =
%f_{0:i\leftarrow 1} \, f_{G:i\leftarrow 0}   =
%f_{0:i\leftarrow 1} u_G \,.
%\]
%We trivially have $f_{G:i\leftarrow 0}  = f_{0:i\leftarrow 0} u_G$
%since $f_{0:i\leftarrow 0}=1$.  So, indeed,
%\[
%f_{G:i\leftarrow b}  =
%f_{0:i\leftarrow b} \, u_G
%\]
%for each $b\in\{0,1\}$.
%
%If (\ref{eq:degeneracy-zeros}) holds for $G$, then $u_G=0$, so again
%\[
%f_{G:i\leftarrow 1}  =   f_{0:i\leftarrow 1} u_G ,
%\]
%and along with the trivial $b=0$ case, we again have (\ref{eq:reduce-f-to-u-soln}) for
%each $b$.  So in fact (\ref{eq:reduce-f-to-u-soln}) holds for each $b$ and all $G$.
%By Theorem \ref{thm:reduce-f-to-u-soln}, this establishes
%degeneracy.
%\eopf   \\

We see from these results
that \textit{any two} distinct minor operations give the \textit{same} notion of
degeneracy, and hence the same goes for \textit{any set} of minor operations.

It is instructive to consider the case of binary matroids,
with $f$ now the indicator function of the cocircuit space.
A loop is an element $e_i$ that belongs to no member of
the cocircuit space.  So, $f(X)=0$ if $e_i\in X$; expressed in the above manner, this
is $f_{G:i\leftarrow 1}=0$ for all $G$, including the case $G=0$ where
$f_{0:i\leftarrow 1}=0$ (i.e., a loop itself is not a member of the cocircuit space).
This ensures that (\ref{eq:degeneracy-products}) holds.
A coloop is an element $e_i$ such that, for all $X\subseteq E\setminus\{e_i\}$,
$X$ belongs to the cocircuit space if and only if $X\cup\{e_i\}$ does.  This includes the
case $X=\emptyset$, which always belongs to the cocircuit space, so $\{e_i\}$ does too.
It follows that $f_{G:i\leftarrow 1}=f_{G:i\leftarrow 0}$, for all $G$, which includes
$f_{0:i\leftarrow 1}=f_{0:i\leftarrow 0}=1$.
If $f_{G:i\leftarrow 1}=f_{G:i\leftarrow 0}=0$, then (\ref{eq:degeneracy-products}) holds
since both sides are 0.  If $f_{G:i\leftarrow 1}=f_{G:i\leftarrow 0}=1$ then
(\ref{eq:degeneracy-products}) holds since every quantity is 1.
Conversely, suppose (\ref{eq:degeneracy-products}) holds for all $G$.
If $f_{0:i\leftarrow 1}=0$, then (\ref{eq:degeneracy-products}) implies that
for all $G$ we have $f_{G:i\leftarrow 1}=0$, so $e_i$ is a loop.
If $f_{0:i\leftarrow 1}=1$, then (\ref{eq:degeneracy-products}) implies
$f_{G:i\leftarrow 1}=f_{G:i\leftarrow 0}$ for all $G$, so $e_i$ is a coloop.

%This ensures that the ratio property (\ref{eq:degeneracy-ratios}) holds.
%A coloop is an element $e_i$ such that, for all $X\subseteq E\setminus\{e_i\}$,
%$X$ belongs to the cocircuit space if and only if $X\cup\{e_i\}$ does.  This includes the
%case $X=\emptyset$, which always belongs to the cocircuit space, so $\{e_i\}$ does too.
%It follows that $f_{G:i\leftarrow 1}=f_{G:i\leftarrow 0}$, for all $G$.  Either they are both 0,
%satisfying (\ref{eq:degeneracy-zeros}), or they are both 1, satisfying
%(\ref{eq:degeneracy-ratios}).  Conversely, suppose either (\ref{eq:degeneracy-ratios})
%or (\ref{eq:degeneracy-zeros}) holds for all $G$.  If $f_{0:i\leftarrow 1}=0$, then
%for all $G$ we have $f_{G:i\leftarrow 1}=0$ by either (\ref{eq:degeneracy-ratios})
%or (\ref{eq:degeneracy-zeros}), so $e_i$ is a loop.  If $f_{0:i\leftarrow 1}=1$, then
%the ratio in (\ref{eq:degeneracy-ratios}) is 1, so (\ref{eq:degeneracy-ratios})
%and (\ref{eq:degeneracy-zeros}) together stipulate that
%$f_{G:i\leftarrow 1}=f_{G:i\leftarrow 0}$ for all $G$, so $e_i$ is a coloop.

\section{Strict binary representations}
\label{sec:strict-bin-reprns}

We now define our notion of faithful representation, and then
determine when it is possible.   \\

\noindent\textbf{Definition}

A \textit{strict binary representation} of a minor-closed set
$\mathcal{A}$ of alternating dimaps is a triple $(F,\varepsilon,\nu)$
such that
\begin{itemize}
\item[(a)]
$F : \mathcal{A} \rightarrow \{\hbox{binary functions}\}$
\item[(b)]
$\varepsilon=(\varepsilon_G\mid G\in\mathcal{A})$ is a family
of bijections $\varepsilon_G:E(G)\rightarrow E(F(G))$;
\item[(c)]
$\nu\in\mathbb{C}$ with $|\nu|=1$;
%  \item[(d)]    %  Not needed: covered by (b).
%   $\varepsilon_G(e)\in E(F(G))$
%   for all $G\in\mathcal{A}$ and all $e\in E(G)$;
\item[(d)]
$F(G^{\omega}) \simeq L^{[\omega]} F(G)$ for all $G\in\mathcal{A}$;
\item[(e)]
$F(G\mino{\mu}e) \simeq F(G)\minor{[\nu\mu]}\varepsilon_G(e)$ for all
$G\in\mathcal{A}$, $e\in E(G)$ and $\mu\in\{1,\omega,\omega^2\}$.
\end{itemize}

We may interpret this definition as follows.  For each alternating dimap $G$,
its corresponding binary function is $F(G)$, and the correspondence between
edges of $G$ and the elements of the ground set of $F(G)$ is given by $\varepsilon_G$.
We require that triality of alternating dimaps corresponds to the trinity transform
for binary functions (condition (d)) and that minor operations correspond too
(condition (e)).  The role of $\nu$ (condition (c)) is simply to allow us to be a little
relaxed about which binary function reduction is used to represent each alternating
dimap reduction.  There is not much room to move here; $\nu$ captures what little
room to move there is.

\vspace{0.1in}

Let $C_1$ denote the ultraloop.
We write $\mathcal{U}_k=\{iC_1\mid i=0,\ldots,k\}$ and
$\mathcal{U}_{\infty}=\{iC_1\mid i\in\mathbb{N}\cup\{0\}\}$,
where $0C_1$ is the empty alternating dimap.

\begin{thm}
If $\mathcal{A}$ is a minor-closed class of alternating dimaps
which has a strict binary representation then $\mathcal{A}=\emptyset$,
or $\mathcal{A}=\mathcal{U}_k$
for some $k$, or $\mathcal{A}=\mathcal{U}_{\infty}$.
\end{thm}

\pf
Suppose $(F,\varepsilon,\nu)$ is
a strict binary representation of $\mathcal{A}$.

The theorem is immediately true if $\mathcal{A}=\emptyset$.
So suppose $\mathcal{A}\not=\emptyset$.

If $|\mathcal{A}|\ge1$ then, since it is minor-closed, it must
contain the empty alternating dimap $C_0$, and the image $F(C_0)$,
representing $C_0$ as a binary function, must be the binary function
$f:2^{\emptyset}\rightarrow\mathbb{C}$ defined by $f(\emptyset)=1$.

So, if $|\mathcal{A}|=1$ then $\mathcal{A}=\mathcal{U}_0$, and the previous paragraph gives
a strict binary representation of $\mathcal{A}$.

Similarly, if $|\mathcal{A}|\ge2$, then it must contain the
ultraloop $C_1$, since that is the only alternating dimap on one edge.   \\

Claim 1:  The image $F(C_1)$ of the ultraloop $C_1$ is given by
\[
F(C_1)  =  \vr{1}{\sqrt{2}-1} .
\]

Proof:

$F(C_1)$ must be some binary
function $f$ on a singleton ground set, $E=\{e\}$ say, with
$f(\emptyset)=1$ and $f(\{e\})=u$ for some $u\in\mathbb{C}$.
Since $C_1$ is self trial, so must $f$ be (by (d) above).  This means
that its vector form $\mathbf{f}=\textvr{1}{u}$ must
%  satisfy
%  \[
%  M(\omega) \mathbf{f} = \mathbf{f} ,
%  \]
%  i.e., it must
be an eigenvector for eigenvalue 1
of the matrix $M(\omega)$.  Now this matrix has eigenvalues 1 and $\omega$,
and the eigenvectors for the former are the scalar multiples of
\[
\vr{1}{\sqrt{2}-1} .
\]
So this is $F(C_1)$, and $u=\sqrt{2}-1$.  So Claim 1 is proved.   \\

If $|\mathcal{A}|=2$ then $\mathcal{A}$ consists of just the empty alternating
dimap and the ultraloop.  The $F$ given by Claim 1, together with
appropriate identity maps $\varepsilon$ (and, in fact, any $\nu$),
gives a strict binary representation.  So we are done in this case.   \\

It remains to deal with $|\mathcal{A}|\ge3$, when $\mathcal{A}$ contains
at least one alternating dimap on two edges.   \\

Claim 2:
The only binary function $\mathbf{f}$ with the property that 
every reduction, on any of the elements of its ground set,
gives $\mathbf{u}=F(C_1)^{\otimes k}$, is 
$\mathbf{f}=F(C_1)^{\otimes(k+1)}$.   \\

Proof:

Observe that, by Claim 1, $\mathbf{u}=(u_G\mid G\in\{0,1\}^E)$ where
$u_G=(\sqrt{2}-1)^{|G|}$, where $|E|=k$ and $|G|$ is the number of 1s in $G$.

Applying Theorem \ref{thm:reduce-f-to-u-soln},
for all $i\in\{0,\ldots,k\}$,
to $\mathbf{u} = F(C_1)^{\otimes k}$ gives
\[
f_{G:i\leftarrow b}  =
f_{0:i\leftarrow b} \, u_{G} \,.
\]
Hence, for each $i$ and each $G$,
\[
f_{G:i\leftarrow 0}  =
f_{0_{k+1}} \, u_{G}  =
u_{G}  =  (\sqrt{2}-1)^{|G|} .
\]
Now consider $f_{G:i\leftarrow 1}$.  Put $j:=0$ if $i\not=0$ and $j:=1$
otherwise (so $j\not=i$).  Then
\begin{eqnarray*}
f_{G:i\leftarrow 1}  & = &
f_{0_{k}:i\leftarrow 1} \, u_{G}  =
f_{(0_{k-1}:i\leftarrow 1):j\leftarrow 0} \, u_{G}  =
f_{0_{k}:j\leftarrow 0} \, u_{0_{k-1}:i\leftarrow 1} \, u_{G}  =
f_{0_{k+1}} \, u_{0_{k-1}:i\leftarrow 1} \, u_{G}  \\
&  = &
1 \cdot  (\sqrt{2}-1) \cdot (\sqrt{2}-1)^{|G|}  =
(\sqrt{2}-1)^{|G|+1} .
\end{eqnarray*}
It follows that, for all $G'\in\{0,1\}^{k+1}$,
\[
f_{G'} = (\sqrt{2}-1)^{|G'|} .
\]
Therefore 
\[
\mathbf{f} = F(C_1)^{\otimes(k+1)} ,
\]
proving the Claim.   \\

Claim 3:
If $k\ge2$ and every reduction of $G$ is $kC_1$, then $G=(k+1)C_1$.   \\

Before proving the claim, consider the case $k=1$, which it does not cover.
Then every alternating dimap on two edges (of which there are four) has
the claimed property.  Of these, the only self-trial one is $2C_1$.   \\

Proof:

Suppose $k=2$.  If $G$ is connected, then there must be some
$e\in E(G)$ and some $\mu\in\{1,\omega,\omega^2\}$ such that
$G\mino{\mu}e=2C_1$ and is therefore disconnected.  The only way in which
$\mu$-reducing a single edge can disconnect a connected alternating dimap is if the edge
is a proper $\mu^{-1}$-semiloop.  It is easily determined that the
only alternating dimaps on three edges which have this property are those
consisting of two triloops and a semiloop.  These do not have three proper semiloops.
So, although they have the specified property for
\textit{one} of their edges, they do not have it for \textit{all} of their
edges.  So $G$ must be disconnected.  Since $G$ has only three edges, some
component of $G$ must be an ultraloop.  But this disappears when reduced,
so the rest of $G$ must be $2C_1$, so $G=3C_1$.

Now suppose $k\ge3$.  It is impossible for $G$ to be connected, because
no reduction of any edge of any connected alternating dimap can possibly
break it up into three or more components.  So consider the components of $G$.
If any of these is not an ultraloop,
then it has at least two edges, and also is left unchanged by reduction of
any edge in any other component (of which there must be at least one), so we
would have a reduction of $G$ that does not give $kC_1$,
which is a contradiction.  So every component of $G$ must be an ultraloop.
Each of these just disappears on reduction, giving $kC_1$, as desired.

So Claim 3 is proved.   \\

Claim 4:
For all $k\ge0$, either $\mathcal{A}$ has no
members with $k$ edges, or it has just one such member which is $kC_1$.   \\

Proof:

We prove the claim by induction on $k$.

We have seen that this is true already for $k\le1$.

Suppose $k=2$.
Every alternating dimap $G_2$
on two edges has the property that every
reduction of it gives the ultraloop.  Therefore, if $G_2\in\mathcal{A}$
then $F(G_2)=F(C_1)^{\otimes2}$, using Claim 2.  But $F(C_1)^{\otimes2}$ is self-trial,
since $F(C_1)$ is.  Therefore $G_2$ must be self-trial too.  But the only
self-trial alternating dimap on two edges is $2C_1$.
So the only member of $\mathcal{A}$ with two edges
is $2C_1$.

Now suppose it is true regarding members of $\mathcal{A}$ with $k-1$
edges, where $k\ge3$.  We show that it is true for $k$ edges.

If $\mathcal{A}$ has no members with $k-1$ edges, then it can have
no members with $k$ edges either, since it is minor-closed.

If $\mathcal{A}$ has at least one member with $k-1$ edges, then by the
inductive hypothesis it can have only one such member, and
this must be $(k-1)C_1$.  We must show that, if $\mathcal{A}$ has at
least one member with $k$ edges, then it can have only one, and it is $kC_1$.

Let $G$ be a member of $\mathcal{A}$ with $k$ edges.  Since $\mathcal{A}$
is minor-closed and has $(k-1)C_1$ as its only member with $k-1$ edges,
all reductions of $G$ must give $(k-1)C_1$.  So, by the requirements of
a strict binary representation, all reductions of $F(G)$ must give
$F(C_1)^{\otimes(k-1)}$.  This implies that $F(G)=F(C_1)^{\otimes k}$,
by Claim 3.  This completes the proof of Claim 4.   \\

It follows from Claim 4 that $\mathcal{A}$ can only be one of the classes
given in the statement of the theorem.  It remains to establish that
a strict binary representation is possible for each of those classes.
This is routine, using
\[
F(kC_1)  =  \vr{1}{\sqrt{2}-1}^{\otimes k}
\]
for every $k$ for which $kC_1\in\mathcal{A}$.
Let $\varepsilon$ consist just of identity maps.
To show that this does indeed enable a strict binary representation,
use Claims 1--3.  The details are a routine exercise.
\eopf

\section{Future work}

It is possible to develop broader definitions of binary
representations of classes of alternating dimaps.
For example, we could allow the edges of $G$ to be
represented by disjoint subsets of the ground set of $F(G)$
instead of just by distinct single elements.  This suggests the problem of
characterising those minor-closed classes of alternating dimaps that
have binary function representations of a more general type, such
as that suggested above.

There may be other ways of representing alternating dimap reductions in a linear way,
using matrices that do not necessarily commute.  We have not ruled
out the possibility that alternating dimaps may indeed be representable as binary functions,
but using other ways of representing the alternating dimap reductions $\bullet[\mu]e$
(rather than the binary function reductions $\bullet\minor{[\mu]}e$ we have studied so far).

%\vspace{0.3in}
%
%\noindent\textbf{Acknowledgements}   [referee from QJM] ******

\end{document}